\documentclass[12pt]{article}
\usepackage{amssymb,latexsym}
\usepackage{amsmath}
\usepackage{tikz}
\newtheorem{theorem}{Theorem}
\newtheorem{corollary}{Corollary}

\def\beq{\begin{equation}}\def\eeq{\end{equation}}
\def\beqn{\begin{eqnarray}}\def\eeqn{\end{eqnarray}}
\def\pont{\hspace{-6pt}{\bf.\ }}

\def\qed{\ifhmode\unskip\nobreak\fi\quad\ifmmode\Box\else$\Box$\fi}
\def\qed{\ifhmode\unskip\nobreak\hfill$\Box$\medskip\fi
 \ifmmode\eqno{\Box}\fi}

\begin{document}

\title{Red-blue clique partitions and (1-1)-transversals}

\author{Andr\'as Gy\'arf\'as \thanks{Research supported in part by
the OTKA Grant No. K104343.}\\\\[-0.8ex]
\small   Alfr\'ed R\'enyi Institute of Mathematics \\[-0.8ex]\and
Jen\H o Lehel\\\\[-0.8ex]
\small   Alfr\'ed R\'enyi Institute of Mathematics \\[-0.8ex]}

\maketitle

\begin{abstract} Motivated by the problem of Gallai on $(1-1)$-transversals of $2$-intervals, it was proved by the authors in 1969 that if the edges of a complete graph $K$ are colored with red and blue (both colors can appear on an edge) so that there is no  monochromatic induced $C_4$ and $C_5$ then the vertices of $K$ can be partitioned into a red and a blue clique. Aharoni, Berger, Chudnovsky and Ziani recently strengthened this by showing that it is enough to assume that there is no induced monochromatic $C_4$ and there is no induced $C_5$ in {\em one of the colors}. Here this is strengthened further, it is enough to assume that there is no monochromatic induced $C_4$ and there is no $K_5$ on which both color classes induce a $C_5$.

We also answer a question of Kaiser and Rabinovich, giving an example of six $2$-convex sets in the plane such that any three intersect but there is no $(1-1)$-transversal for them.
\end{abstract}

\section{Red-blue clique partitions of complete graphs}

In 1968, thinking on a problem about piercing cycles of digraphs, Gallai arrived to the problem of piercing $2$-intervals.  He defined $2$-intervals as sets of the real line $R$ having two interval components, one in $(-\infty,0)$ and one in $(0,\infty)$ and asked: how many points are needed to pierce a family of pairwise intersecting $2$-intervals? His question generated \cite{GYL1} in which (as a special case of a general upper bound) we proved that two points always pierce pairwise intersecting $2$-intervals and one of them can be selected from $(-\infty,0)$ and the other from $(0,\infty)$. Let's call such a pair of points a $(1-1)$-transversal. This result can be extended to $2$-trees, where a  $2$-tree is the union of two subtrees, one is a subtree of $T_1$ the other is a subtree of $T_2$, where $T_1$ and $T_2$ are  vertex-disjoint trees. 
In \cite{GYL1} we proved a stronger result, Theorem \ref{no45}, using only properties of the intersection graph of subtrees of a tree. Consider {\em $2$-colored complete graphs}, where edges are colored with red, blue, or both colors. Edges of one color only are called {\em pure} edges, they can be pure red or pure blue. Another view is to consider a complete graph (clique) as the union of a red and a blue graph on the same vertex set.

\begin{theorem}\pont(Gy\'arf\'as, Lehel \cite{GYL1}, 1970)\label{no45} Assume that $G$ is a $2$-colored complete graph containing no monochromatic induced $C_4$ and $C_5$. Then $V(G)$ can be partitioned into a red and a blue clique.
\end{theorem}

Given a set of $n$ pairwise intersecting $2$-subtrees, one can represent their intersections by a $2$-colored complete graph $K_n$. Then both colors determine {\em chordal graphs} i.e. graphs in which every cycle of length at least four has a chord.
In particular, there is no monochromatic induced $C_4$ or $C_5$. Applying Theorem \ref{no45}, the vertices of $K_n$ can be partitioned into a red and a blue clique (empty sets or one vertex is accepted as a clique) and by the Helly-property of subtrees we have a $(1-1)$-transversal for the $2$-subtrees. Thus Theorem \ref{no45} implies the following.

\begin{corollary}\pont(Gy\'arf\'as, Lehel \cite{GYL1}, 1970)\label{tree} Pairwise intersecting $2$-subtrees have a $(1-1)$ transversal.
\end{corollary}

Since $2$-colorings of complete graphs with pure edges only can be considered as a graph and its complement, we get another consequence of Theorem \ref{no45}.

\begin{corollary}\pont (F\"oldes, Hammer \cite{FH}, 1977)\label{fh} Assume that a graph $G$ does not contain $C_4,2K_2,C_5$ as an induced subgraph. Then $G$ is a split graph, i.e. its vertices can be partitioned into a clique and an independent set.
\end{corollary}

The seminal paper of Tardos  \cite{TA} (1995) introduced topological methods, he proved that $2$-intervals without $k+1$ pairwise disjoint members have $(k-k)$-transversals. Methods of Kaiser \cite{KA} (1997), Alon \cite{A1,A2} (1998, 2002), Matousek \cite{MAT} (2001), Berger \cite{BER} (2005)  brought many nice results and this list of references is very far from being complete.  In this note we only consider the graph coloring approach.
Very recently Theorem \ref{no45} was generalized as follows.

\begin{theorem}\pont (Aharoni, Berger, Chudnovsky, Ziani \cite{ABCZ}, 2015) \label{be} Assume that $G$ is a $2$-colored complete graph such that there is no monochromatic induced $C_4$ and there is no red induced $C_5$.  Then $V(G)$ can be partitioned into a red and a blue clique.
\end{theorem}

We show that the proof of Theorem \ref{no45} in \cite{GYL1} yields an even stronger result. Let $K_5^*$ denote the $2$-colored $K_5$ where every edge is pure and both colors span a $C_5$.

\begin{theorem}\pont  \label{c4plus} Assume that $G$ is a $2$-colored complete graph such that there is no monochromatic induced $C_4$ and there is no $K^*_5$. Then $V(G)$ can be partitioned into a red and a blue clique.
\end{theorem}

\noindent
\bf Proof. \rm  We prove by induction on $|V(G)|$, for $1\le |V(G)|\le 3$ the theorem is obvious. Fixing any $p\in V(G)$, by the inductive hypothesis we have $V(G-p)=R\cup B$ where $R$ and $B$ are disjoint vertex sets spanning a red and a blue clique.

Set $$R^*=\{r\in R: (p,r) \mbox{ is pure blue} \}, B^*=\{b\in B: (p,b) \mbox{ is pure red} \}.$$
Assume that among all choices of $R,B$ satisfying $V(G-p)=R\cup B$, $|R^*|+|B^*|$ is as small as possible.
We show that either $R^*$ or $B^*$ is empty, thus $R$ or $B$ can be extended with $p$, concluding the proof.

Suppose on the contrary that $R^*,B^*$ are both nonempty.  For any  $q\in B^*$ there exists $r\in R$ such that $(q,r)$ is pure blue, otherwise $R_1=R\cup \{q\}$ and $B_1=B\setminus \{q\}$ would be a red-blue clique partition of $V(G-p)$ with  $|R_1^*|+|B_1^*|<|R^*|+|B^*|$, contradicting the assumption. In fact, we may assume that $r\in R^*$, otherwise, with any $s\in R^*$, consider the four-cycle  $C=(p,q,s,r,p)$. If $(q,s)$ would be red then $C$ is a red cycle with pure blue diagonals $(q,r),(p,s)$, contradiction. Thus $(q,s)$ is pure blue and we can choose $s\in R^*$ instead of $r$.

Applying the argument of the previous paragraph for any $s\in R^*$, there exists $q\in B^*$ such that $(s,q)$ is pure red. Thus there exists a shortest even cycle $C=(s_1,q_1,s_2,\dots,q_m,s_1)$ in the bipartite graph $[R^*,B^*]$ with edges alternating as pure red, pure blue, pure red... We claim that $C$ is a four-cycle. Indeed, if $m>2$, then from the minimality of $m$, all diagonals  $(s_i,q_j)$ must have both colors. In particular, $(s_1,q_2),(s_3,q_1)$ both have red colors. Now if $(q_1,q_2)$ is pure blue then the red four-cycle $(s_1,q_1,p,q_2,s_1)$ has pure blue diagonals, otherwise the red four-cycle $(q_1,q_2,s_2,s_3,q_1)$ has pure blue diagonals, contradicting the assumption that there is no induced monochromatic $C_4$. This proves the claim, $C=(s_1,q_1,s_2,q_2,s_1)$. Observe that $(s_1,s_2)$ is pure red, otherwise $(s_1,s_2,q_1,q_2,s_1)$ is a blue four-cycle with pure red diagonals and $(q_1,q_2)$ is pure blue, otherwise $(q_1,q_2,s_2,s_1,q_1)$ is a red four-cycle with pure blue diagonals, contradiction.

Therefore $\{p,s_1,s_2,q_1,q_2\}$ spans a $K_5^*$, giving the final contradiction. \qed

\section{Red-blue clique partition of complete $3$-uniform hypergraphs}

Extending $2$-intervals Kaiser and Rabinovich [10]
defined a planar $2$-body as a union of two closed convex sets of the plane 
separated by a fixed line, say the $y$-axis.  They asked whether 
the assumption `any three $2$-bodies intersect' implies that they have a $(1-1)$-transversal.

Following Theorem \ref{no45}, where the $(1-1)$-transversal is translated
into properties that imply a red-blue clique cover of a two-colored clique, this problem can be stated in terms of red-blue colored $3$-uniform complete hypergraphs. 
However, the obstructions for a red-blue clique cover of a hypergraph can be more complicated than those for graphs in 
Theorems \ref{no45}, \ref{be} and \ref{c4plus}. In particular, as our next example shows, 
the answer is negative to the 
question above.\\

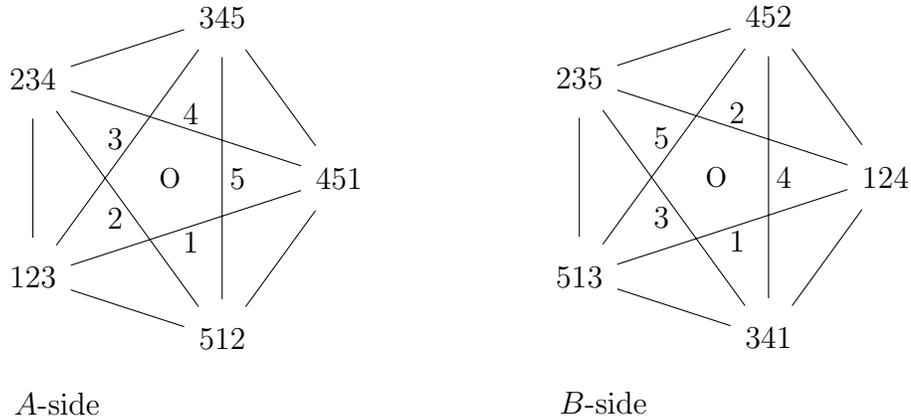
\begin{figure}
\hskip1.5truecm
\begin{tikzpicture}[scale=1.5]
\tikzstyle{every node}=[draw=white,shape=circle];
\tikzstyle{txt}  = [circle, minimum width=1pt, draw=white, inner sep=0pt]
\path (0:0cm) node (v0) {\Large o};
\path (0:.6cm) node (v5) {$5$};
\path (72:.6cm) node (v4) {$4$};
\path (2*72:.6cm) node (v3) {$3$};
\path (3*72:.6cm) node (v2) {$2$};
\path (4*72:.6cm) node (v1) {$1$};
\path (0:1.5cm) node (p1) {$451$};
\path (72:1.5cm) node (p2) {$345$};
\path (2*72:1.5cm) node (p3) {$234$};
\path (3*72:1.5cm) node (p4) {$123$};
\path (4*72:1.5cm) node (p5) {$512$};
\draw (p1) -- (p2) -- (p3) -- (p4) -- (p5) -- (p1);
\draw (p1) -- (p3) -- (p5) -- (p2) -- (p4) -- (p1);
\node[txt] () at (-1,-2){$A$-side};
\end{tikzpicture}
\hskip2truecm
\begin{tikzpicture}[scale=1.5]
\tikzstyle{every node}=[draw=white,shape=circle];
\tikzstyle{txt}  = [circle, minimum width=1pt, draw=white, inner sep=0pt]
\path (0:0cm) node (v0) {\Large o};
\path (0:.6cm) node (v4) {$4$};
\path (72:.6cm) node (v2) {$2$};
\path (2*72:.6cm) node (v5) {$5$};
\path (3*72:.6cm) node (v3) {$3$};
\path (4*72:.6cm) node (v1) {$1$};
\path (0:1.5cm) node (p1) {$124$};
\path (72:1.5cm) node (p2) {$452$};
\path (2*72:1.5cm) node (p3) {$235$};
\path (3*72:1.5cm) node (p4) {$513$};
\path (4*72:1.5cm) node (p5) {$341$};
\draw (p1) -- (p2) -- (p3) -- (p4) -- (p5) -- (p1);
\draw (p1) -- (p3) -- (p5) -- (p2) -- (p4) -- (p1);
\node[txt] () at (-1,-2){$B$-side};
\end{tikzpicture}
\caption{Any three of the $2$-bodies, $A_i\cup B_i$, $0\leq i\leq 5$, intersect,
and they do not have a $(1-1)$-transversal}
\end{figure}

\noindent {\bf Example 1.} We define six planar $2$-bodies, $A_i\cup B_i$, $0\leq i\leq 5$, as follows. On each side of the $y$-axis we are given $5$ triangles formed by consecutive triples of vertices of a fixed regular pentagon, and the inner pentagon bordered by its diagonals is the sixth convex set. On the $A$-side  the (clockwise) consecutive triangles  are labeled 
$A_1, A_2, A_3, A_4, A_5$; on the $B$-side 
the labeling of the consecutive triangles  is $B_1, B_3,B_5, B_2, B_4$; the inner pentagons are labeled $A_0$ and $B_0$ (see Fig.1).\\

The $2$-bodies of the example define a natural $2$-coloring of the edges of $K_6^{(3)}$, the complete $3$-uniform hypergraph on vertex set $V=\{0,1,2\dots,5\}$:  if a triple of convex sets has
non-empty intersection on the $A$-side (on the $B$-side), then the corresponding edge
of the hypergraph is colored red  (blue). It is easy to check that  no four convex sets intersect on either side, furthermore, the $10$ vertices and the $10$ intersection points of diagonals are the intersections of the triples of the six $2$-bodies.
The red edges are the triples $\{(i,i+1,i+2): 1\le i \le 5\}$ (counting modulo $5$) and their complements (with respect to $V$); the blue edges are the triples  $\{(i,i+1,i+3 ): 1\le i \le 5\}$ and their complements. Thus $V$ is covered by two triples of the same color, but it cannot be covered by a red and a blue triple.\\

 Example 1 shows that the assumption `any three $2$-bodies intersect'
does not imply that the $2$-bodies have a $(1-1)$-transversal. However, Kaiser and Rabinovich \cite{KARA}  proved it from the condition that `any four $2$-bodies intersect'.

\begin{theorem}\pont (Kaiser, Rabinovich \cite{KARA}, 1999) \label{34} Assume that $S$ is a set of planar $2$-bodies such that any four members of $S$ have nonempty intersection. Then $S$ has a $(1-1)$-transversal.
\end{theorem}

\end{document}